\documentclass[12pt]{amsart}

\usepackage{float}
\usepackage{geometry}
\usepackage[all]{xy}
\usepackage{amssymb, amsmath}
\usepackage[breaklinks]{hyperref}
\usepackage{mathtools}
\usepackage{xr}
\DeclareMathAlphabet{\mathpzc}{OT1}{pzc}{m}{it}

\newtheorem{thm}{Theorem}[section]

\newtheorem{prob}{Problem}

\theoremstyle{definition}

\theoremstyle{remark}

\numberwithin{figure}{section}%  Fig. 4.1
\numberwithin{table}{section}
\numberwithin{equation}{section}
\newcommand{\eps}{\varepsilon}
\newcommand{\To}{\longrightarrow}

\newcommand{\restn}{\mathrm{Rest}_{x_n=0}}

\mathchardef\mhyphen="2D

           \newcommand{\eq}[1][r]
   {\ar@<-3pt>@{-}[#1]
    \ar@<-1pt>@{}[#1]|<{}="gauche"
    \ar@<+0pt>@{}[#1]|-{}="milieu"
    \ar@<+1pt>@{}[#1]|>{}="droite"
    \ar@/^2pt/@{-}"gauche";"milieu"
    \ar@/_2pt/@{-}"milieu";"droite"}

\newcommand{\iotan}{\iota_{\frac{\partial}{\partial x_n}}}

\newcommand{\R}{\mathbb R}  %@@
\newcommand{\Z}{\mathbb Z}  %@@
\newcommand{\N}{\mathbb N}  %@@
\newcommand{\C}{\mathbb C}  %@@
\newcommand{\rightsetse}[1]{%
\hidewidth\rotatebox[origin=c]{-45}{$\xrightarrow{\kern2em}$}
     \rlap{\raisebox{1ex}
     {$\kern-.8em\scriptstyle #1$}}\hidewidth}
\newcommand{\rightsetsw}[1]{%
\hidewidth\rotatebox[origin=c]{45}{$\xleftarrow{\kern2em}$}
     \rlap{\raisebox{.1ex}
     {$\kern-.8em\scriptstyle #1$}}\hidewidth}
\newcommand{\leftsetsw}[1]{%
\hidewidth
     \llap{\raisebox{1ex}
     {$\scriptstyle #1$\kern-.8em}}
    \rotatebox[origin=c]{45}{$\xleftarrow{\kern2em}$}\hidewidth}
\newcommand{\rightsetnw}[1]{%
\hidewidth\rotatebox[origin=c]{135}{$\xrightarrow{\kern2em}$}
     \rlap{\raisebox{1ex}
     {$\kern-.8em\scriptstyle #1$}}\hidewidth}
\newcommand{\rightsetd}[1]{%
\hidewidth\rotatebox[origin=c]{-90}{$\xrightarrow{\kern2em}$}
     \rlap{{$\scriptstyle #1$}}\hidewidth}
\usepackage[breaklinks]{hyperref}
\begin{document}

\title[Classification of Symmetry breaking operators on spheres]
{Classification of differential symmetry breaking operators for differential forms}

\author{Toshiyuki KOBAYASHI, Toshihisa KUBO, Michael PEVZNER}

\thanks{\emph{Acknowledgements}:
The first named author was partially supported by Institut des Hautes \'Etudes Scientifiques, France and Grant-in-Aid for Scientific Research (A) (25247006), Japan Society for the Promotion of
Science. All three authors were partially supported by CNRS Grant PICS n$^\mathrm{o}$ 7270.}

%\date{\today}

%%%%%%%%%%%%%%%%%%%%%%%%%%%%%%%%%%%%%
\begin{abstract}
We give a complete classification of 
 conformally covariant differential operators between the spaces of differential $i$-forms on the sphere $S^n$ and $j$-forms on the totally geodesic hypersphere $S^{n-1}$ by analyzing the restriction of principal series representations of the Lie group $O(n+1,1)$. 
 Further, we provide explicit formul\ae{} for these matrix-valued operators in the flat coordinates and find  factorization identities for them. 

 \vskip5pt
\noindent This note was published in \href{http://dx.doi.org/10.1016/j.crma.2016.04.012}{ C. R. Acad. Sci. Paris, Ser. I, (2016),}  
{http://dx.doi.org/
10.1016/j.crma.2016.04.012}.
\vskip7pt

\end{abstract}

\maketitle

%%%%%%%%%%%%%%%%%%%%%%%%%%%%%%%%%%%%%%%%%%%%%%%
\section{Introduction}
Suppose a Lie group $G$ acts conformally on a Riemannian manifold $(X,g)$. 
This means that there exists a positive-valued function $\Omega\in C^\infty(G\times X)$ (\emph{conformal factor}) such that
$$
L_h^*g_{h\cdot x}=\Omega(h,x)^2g_x \quad\mathrm{for\; all}\;h\in G\;\mathrm{and}\; x\in X,
$$
where $L_h:X\to X, x\mapsto h\cdot x$ denotes the action of $G$ on $X$. 
Since $\Omega$ satisfies a cocycle condition, 
we can form a family of representations $\varpi^{(i)}_{u}$ for $u\in\C$
 and $0\leq i\leq \dim X$ on the space $\mathcal E^i(X)$ of differential $i$-forms on $X$ by
\begin{equation}\label{eqn:varpi}
\varpi^{(i)}_{u}(h)\alpha:=\Omega(h^{-1},\cdot)^u 
L_{h^{-1}}^*\alpha\quad (h\in G).
\end{equation}
The representation $\varpi^{(i)}_{u}$ of the conformal group $G$ on $\mathcal E^i(X)$
will be simply denoted by
 $\mathcal E^i(X)_{u}$.

If $Y$ is a submanifold of $X$, 
then we can also define a family of representations $\varpi^{(j)}_{v}$ on $\mathcal E^j(Y)$ ($v\in\C, 0\leq j\leq\dim Y$) of the subgroup
$$
G':=\left\{h\in G\,:\,h\cdot Y=Y\right\},
$$
which acts conformally on the Riemannian submanifold $(Y,g{\;\vert}_{ Y})$.

We study  
differential operators $\mathcal D:\mathcal E^i(X)\To\mathcal E^j(Y)$ that intertwine
the two representations $\varpi^{(i)}_{u}\vert_{G'}$ and $\varpi^{(j)}_{v}$ of $G'$.
Here $\varpi^{(i)}_{u}\vert_{G'}$ stands for the restriction of the $G$-representation
$\varpi^{(i)}_{u}$ to the subgroup $G'$. We say that such $\mathcal D$ is a 
\emph{differential symmetry breaking operator},
and denote by
$\mathrm{Diff}_{G'}(\mathcal E^i(X)_{u},\mathcal E^j(Y)_{v})$
the space of all differential symmetry breaking operators.
We address the following problems:
\vskip 0.1in

\begin{prob}\label{prob:1}
Determine the dimension
of the space $\mathrm{Diff}_{G'}\left(\mathcal{E}^i(X)_{u}, \mathcal{E}^j(Y)_{v}\right)$.
In particular, find a necessary and sufficient 
condition on a quadruple $(i,j,u,v)$ such that there exist nontrivial
differential symmetry breaking operators. 
\end{prob}

\begin{prob}\label{prob:2}
Construct explicitly a basis of 
$\mathrm{Diff}_{G'}\left(\mathcal{E}^i(X)_{u}, \mathcal{E}^j(Y)_{v}\right)$.
\end{prob}

In the case where $X=Y$, $G=G'$, and $i = j = 0$,
a classical prototype of such operators is 
a second order differential operator called the
Yamabe operator
$$
\Delta+\frac{n-2}{4(n-1)}\kappa\in \mathrm{Diff}_{G}(\mathcal E^0(X)_{\frac n2-1},\mathcal E^0(X)_{\frac n2+1}),
$$
$\Delta$ is the Laplace--Beltrami operator, where $n$ is the dimension of $X$, and $\kappa$ is the
scalar curvature of $X$. 
Conformally covariant differential operators of higher order are also known:
the Paneitz operator (fourth order) \cite{P08}, 
which appears in four dimensional supergravity \cite{FT},
or more generally, the so-called 
GJMS operators \cite{GJMS} are such examples.
Analogous conformally covariant operators on  forms  
($i=j$ case) were studied 
by Branson \cite{Branson}. 
On the other hand, the insight of representation theory of conformal groups
is useful in studying Maxwell's equations, see \cite{KW14}, for instance.

Let us consider the more general case where $Y\neq X$ and $G' \neq G$. 
An obvious example of symmetry breaking operators is the restriction operator
$\mathrm{Rest}_Y$ which belongs to 
$\mathrm{Diff}_{G'}\left(\mathcal{E}^i(X)_{u}, \mathcal{E}^i(Y)_{u}\right)$
for all $u\in\C$. Another elementary example
is $\mathrm{Rest}_Y\circ \iota_{N_Y(X)} \in 
\mathrm{Diff}_{G'}\left(\mathcal{E}^i(X)_{u},
\mathcal{E}^{i-1}(Y)_{v}\right)$
if $v = u+1$ where $\iota_{N_Y(X)}$ denotes
the interior multiplication by the normal vector field to $Y$
when $Y$ is of codimension one in $X$.

In the model space where $(X,Y)=(S^n,S^{n-1})$, the pair $(G,G')$ of conformal groups amounts
to $(O(n+1,1),O(n,1))$ modulo center, and
Problems \ref{prob:1} and \ref{prob:2} have been recently solved for $i=j=0$ 
by Juhl \cite{Juhl}, see also \cite{K14,KP2} and \cite{KOSS15}
for different approaches by the F-method and the 
residue calculus, respectively. 

 Problems \ref{prob:1} and \ref{prob:2}  for general $i$ and $j$ for the model space can be reduced to analogous problems for (nonspherical) principal series representations by the isomorphism \eqref{eqn:wI} below. In this note
 we shall give complete solutions to Problems \ref{prob:1} and \ref{prob:2}  in those terms (see Theorems \ref{thm:A} and \ref{thm:B}).

Notation: $\N=\{0,1,2,\cdots\}$, $\N_+=\{1,2,\cdots\}$.

%%%%%%%%%%%%%%%%%%%%%%%%%%%%%%%%%%%%%%%%%%%%%%%
 \section{Principal series representations of $G=O(n+1,1)$}
 
 We set up notations. Let $P=MAN$ be a Langlands decomposition of a minimal parabolic subgroup of
 $G=O(n+1,1)$. For $0\leq i\leq n$, $\delta\in \Z/2\Z$, and $\lambda\in\C$, we extend the outer tensor product representation $\bigwedge^i(\C^n)\otimes(-1)^\delta\otimes\C_\lambda$ of $MA\simeq
 (O(n)\times O(1))\times\R$ to $P$ by letting $N$ act trivially, and form a $G$-equivariant
 vector bundle 
$
\mathcal V^i_{\lambda,\delta}:=G\times_P
\left(\bigwedge^i(\C^n)\otimes(-1)^\delta\otimes\C_\lambda\right)
$
 over the real flag variety
$X=G/P\simeq S^n$. Then we define an unnormalized principal
 series representations
 \begin{equation}\label{eqn:psr}
 I(i,\lambda)_\delta:=\mathrm{Ind}_P^G\left(
\bigwedge^i(\C^n)\otimes(-1)^\delta\otimes\C_\lambda
\right)
 \end{equation}
of $G$ on the Fr\'echet space $C^\infty(X,\mathcal V_{\lambda,\delta}^i)$ of smooth sections.

In our parametrization, $I(i,n-2i)_\delta$ and $I(i,i)_\delta$ have the same infinitesimal character with the trivial one-dimensional representation of $G$.
Then, for all $u\in\C$, we have a natural $G$-isomorphism
\begin{equation}\label{eqn:wI}
\varpi^{(i)}_u\simeq I(i,u+i)_{i\,\mathrm{mod}\,2}.
\end{equation}
Similarly, for $0\leq j\leq n-1$, $\eps\in \Z/2\Z$ and $\nu\in\C$, we define an unnormalized principal series representation
$
 J(j,\nu)_\eps:=\mathrm{Ind}_{P'}^{G'}\left(
\bigwedge^j(\C^{n-1})\otimes(-1)^\eps\otimes\C_\nu
\right)$
 of the subgroup $G'=O(n,1)$  
on $C^\infty(Y,\mathcal W_{\nu,\eps}^j)$, where $\mathcal W_{\nu,\eps}^j:=G'\times_{P'}
\left(\bigwedge^j(\C^{n-1})\otimes(-1)^\eps\otimes\C_\nu\right)$  is a $G'$-equivariant vector bundle 
over $Y=G'/P'\simeq S^{n-1}$.

%%%%%%%%%%%%%%%%%%%%%%%%%%%%%%%%%%%%%%%%%%%%%%%
\section{Existence condition for differential symmetry breaking operators}

A continuous $G'$-intertwining operator $T:I(i,\lambda)_\delta\To J(j,\nu)_\eps$ is said to be a 
\emph{symmetry breaking operator} (SBO). We say that $T$ is a differential operator if $T$ satisfies
$\mathrm{Supp}(Tf)\subset \mathrm{Supp} f$ for all $f\in C^\infty(X,\mathcal V^i_{\lambda,\delta})$,
and $\mathrm{Diff}_{G'}\left( I(i,\lambda)_\delta, J(j,\nu)_\eps\right)$ denotes the space of differential
SBOs. We give a complete solution to Problem \ref{prob:1} for $(X,Y)=(S^n,S^{n-1})$ in terms of principal
series representations:

\begin{thm}\label{thm:A}
Let $n\geq 3$. Suppose $0\leq i\leq n$, $0\leq j\leq n-1$, $\lambda,\nu\in\C$,
and $\delta,\eps\in\Z/2\Z$. 
 Then the following
three conditions on 6-tuple $(i,j,\lambda,\nu,\delta,\eps)$ are equivalent:
\begin{enumerate}
\item[\emph{(i)}]
$\mathrm{Diff}_{O(n,1)}(
I(i,\lambda)_\delta,
J(j,\nu)_{\eps})\neq\{0\}.$
\item[\emph{(ii)}]
$\dim\mathrm{Diff}_{O(n,1)}(
I(i,\lambda)_\delta,
J(j,\nu)_{\eps})=1.$
\item[\emph{(iii)}]
The 6-tuple belongs to one of the following six cases:
\vskip 0.05in
\begin{enumerate}
\item[] \emph{Case 1}. $j = i$, 
$0\leq i \leq n-1$, $\nu - \lambda \in \N$, 
$\eps-\delta \equiv \nu -\lambda \; \mathrm{mod}\, 2$.
\vskip 0.05in

\item[] \emph{Case 2}. $j = i-1$, 
$1 \leq i \leq n$, $\nu-\lambda \in \N$, 
$\eps -\delta \equiv \nu - \lambda \; \mathrm{mod}\, 2$.
\vskip 0.05in

\item[] \emph{Case 3}. $j=i+1$, 
$1 \leq i \leq n-2$, 
$(\lambda,\nu) = (i,i+1)$, 
$\eps \equiv \delta + 1 \; \mathrm{mod}\, 2$.
\vskip 0.05in

\item[] \emph{Case 3$'$}. $(i,j) = (0,1)$, 
$-\lambda \in \N$, $\nu = 1$, $\eps \equiv \delta + \lambda+1 \; \mathrm{mod}\, 2$.
\vskip 0.05in

\item[] \emph{Case 4}. $j=i-2$,
$2\leq i \leq n-1$, $(\lambda,\nu) = (n-i,n-i+1)$,
$\eps  \equiv \delta + 1\; \mathrm{mod}\, 2$.
\vskip 0.05in

\item[] \emph{Case 4$'$}. $(i,j) = (n, n-2)$,
$-\lambda \in \N$, $\nu=1$, $\eps \equiv \delta + \lambda+1\; \mathrm{mod}\, 2$.
\vskip 0.05in
\end{enumerate}
\end{enumerate}
\end{thm}
\noindent
We set $\Xi:=\{(i,j,\lambda,\nu)$: the 6-tuple $(i,j,\lambda,\nu,\delta,\eps)$ satisfies one of the equivalent conditions of Theorem \ref{thm:A} for some $\delta,\eps\in\Z/2\Z\}$.

%%%%%%%%%%%%%%%%%%%%%%%%%%%%%%%%%%%%%%%%%%%%%%%
\section{Construction of differential symmetry breaking operators}
In this section, we describe an explicit generator of the space of differential SBOs if one of the equivalent conditions
in Theorem \ref{thm:A} is satisfied. For this we use the \emph{flat picture} of the principal
series representations $I(i,\lambda)_\delta$ of $G$
which realizes the representation space $C^\infty(X,\mathcal V^i_{\lambda,\delta})$ as a subspace of
 $C^\infty(\R^n,\bigwedge^i(\C^n))$ by 
trivializing the bundle $\mathcal V^i_{\lambda,\delta}\To X$ on the open Bruhat cell
$$
\R^n\hookrightarrow X,\quad (x_1,\cdots,x_n)\mapsto\exp\left(\sum_{j=1}^nx_jN_j^-\right)P.
$$
Here $\{N_1^-,\cdots,N_n^-\}$ is an orthonormal basis of 
the nilradical $\mathfrak{n}_-(\R)$ of the opposite parabolic
subalgebra with respect to an $M$-invariant inner product.
Without loss of generality, we may and do assume that the open Bruhat cell $\R^{n-1}\hookrightarrow Y\simeq G'/P'$
is given by putting $x_n=0$. Then the flat picture of the principal series representation $J(j,\nu)_\eps$
of $G'$ is defined by realizing 
$C^\infty(Y,\mathcal W_{\nu,\eps}^j)$ as a subspace of 
$C^\infty(\R^{n-1},\bigwedge^j(\C^{n-1}))$.
For the construction of explicit generators of matrix-valued SBOs,
we begin with a scalar-valued differential operator.
 For $\alpha\in\C$ and $\ell\in\N$, we define a polynomial of two variables $(s,t)$ by
$$
\left( I_\ell\widetilde C^\alpha_\ell\right)(s,t):=s^{\frac\ell 2}\widetilde C^\alpha_\ell\left(\frac{t}{\sqrt s}\right),
$$
where $\widetilde C^\alpha_\ell(z)$ is the renormalized Gegenbauer polynomial given by
$$
\widetilde C^\alpha_\ell(z):=\frac1{\Gamma\left(\alpha+\left[\frac{\ell+1}2\right]\right)}
\sum_{k=0}^{\left[\frac{\ell}2\right]}(-1)^k
\frac{\Gamma (\ell-k+\alpha)}{k!(\ell-2k)!}(2z)^{\ell-2k}.
$$
Then $\widetilde C^\alpha_\ell(z)$ is a nonzero polynomial for all $\alpha\in\C$ and $\ell\in\N$, and
a (normalized) \emph{Juhl's conformally covariant operator} $\widetilde \C_{\lambda,\nu}: C^\infty(\R^n)\To C^\infty(\R^{n-1})$
is defined by
$$
\widetilde \C_{\lambda,\nu}:=\restn\circ\left( I_\ell\widetilde C^{\lambda-\frac{n-1}2}_\ell\right)
\left(-\Delta_{\R^{n-1}},\frac\partial{\partial x_n}\right),
$$
for $\lambda,\nu\in\C$ with $\ell:=\nu-\lambda\in\N$.
For instance, 
$$
\widetilde \C_{\lambda,\nu}=\mathrm{Rest}_{x_n=0}\circ
\begin{cases}
\mathrm{id} & \mathrm{if}\;  \nu=\lambda,\\
2\frac{\partial}{\partial x_n} & \mathrm{if}\; \nu=\lambda+1,\\
\Delta_{\R^{n-1}}+(2\lambda-n+3)
\frac{\partial^2}{\partial x_n^2}& \mathrm{if}\; \nu=\lambda+2.
%\\2\Delta_{\R^{n-1}}\frac{\partial}{\partial x_n}+\frac23(2\lambda-n+5)\frac{\partial^3}{\partial x_n^3}
%& \mathrm{if}\; \nu=\lambda+3.
\end{cases}
$$
For $(i,j,\lambda,\nu)\in\Xi$, we introduce a new family of matrix-valued differential operators
$$
\widetilde\C^{i,j}_{\lambda,\nu}: C^\infty(\R^n,\bigwedge^i(\C^n))\To C^\infty(\R^{n-1},\bigwedge^j(\C^{n-1})),
$$
by using the identifications $\mathcal E^i(\R^n)\simeq C^\infty(\R^n)\otimes\bigwedge^i(\C^n)$
and $\mathcal E^j(\R^{n-1})\simeq C^\infty(\R^{n-1})\otimes\bigwedge^j(\C^{n-1})$,
as follows. Let $d^*_{\R^n}$ be the codifferential,
which is the formal adjoint of the differential $d_{\R^n}$,
 and $\iotan$ the inner multiplication by the vector field
$\frac\partial{\partial x_n}$. Both operators map $\mathcal E^i(\R^n)$ to $\mathcal E^{i-1}(\R^n)$. 
For $\alpha\in\C$ and $\ell\in\N$, let
$
\gamma(\alpha,\ell):=
1$ ($\ell$ is odd); $=
\alpha+\frac \ell2$ ($ \ell$ is even).
Then we set
\begin{align}
\C^{i,i}_{\lambda,\nu}
&:=\widetilde{\C}_{\lambda+1, \nu-1} d_{\R^n} d^*_{\R^n}
-\gamma(\lambda-\frac{n}{2}, \nu-\lambda) \widetilde{\C}_{\lambda,\nu-1}
d_{\R^n}\iota_{\frac{\partial}{\partial x_n}} 
+ \frac{1}{2}(\nu-i)\widetilde{\C}_{\lambda,\nu} \nonumber \qquad \qquad \qquad \mathrm{for}\; 0\leq i \leq n-1.
\\
\C^{i,i-1}_{\lambda,\nu} 
&:= - \widetilde{\mathcal{\C}}_{\lambda+1, \nu-1} d_{\R^n}d^*_{\R^n}
\iota_{\frac{\partial}{\partial x_n}} -\gamma(\lambda-\frac{n-1}{2}, \nu-\lambda)
\widetilde{\C}_{\lambda+1, \nu}d^*_{\R^n} 
+\frac{1}{2}(\lambda+i-n)\widetilde{\C}_{\lambda,\nu} \iota_{\frac{\partial}{\partial x_n}}
\; \mathrm{for}\; 1\leq i \leq n.
\nonumber 
\end{align}
We note that there exist isolated parameters 
$(\lambda,\nu)$ for which $\C^{i,i}_{\lambda,\nu} =0$
or $\C^{i,i-1}_{\lambda,\nu} =0$. 
For instance, $\C^{0,0}_{\lambda,\nu} = \frac{1}{2}\nu \widetilde{\C}_{\lambda,\nu}$,
and thus $\C^{0,0}_{\lambda,\nu} = 0$ if $\nu =0$.
To be precise, we have the following:
\begin{itemize}
\item[]$\C^{i,i}_{\lambda, \nu} = 0$ 
if and only if $\lambda = \nu =i$ or $\nu = i = 0$;
\item[]$\C^{i,i-1}_{\lambda,\nu}=0$ if and only if 
$\lambda=\nu=n-i$ or $\nu = n-i = 0$.
\end{itemize}

We renormalize these operators by
\begin{equation*}
\widetilde{\C}^{i,i}_{\lambda,\nu}:=
\begin{cases}
\mathrm{Rest}_{x_n=0} & \text{if $\lambda = \nu$},\\
\widetilde{\C}_{\lambda,\nu} & \text{if $i=0$},\\
\C^{i,i}_{\lambda,\nu} & \text{otherwise},
\end{cases}
\quad
\text{and}
\quad
\widetilde{\C}^{i,i-1}_{\lambda,\nu}:=
\begin{cases}
\mathrm{Rest}_{x_n=0} \circ \iota_{\frac{\partial}{\partial x_n}} & \text{if $\lambda = \nu$},\\
\widetilde{\C}_{\lambda,\nu} \circ \iota_{\frac{\partial}{\partial x_n}} & \text{if $i=n$},\\
\C^{i,i-1}_{\lambda,\nu} & \text{otherwise}.
\end{cases}
\end{equation*}

\noindent
Then $\widetilde{\C}^{i,i}_{\lambda,\nu}$ $(0\leq i \leq n-1)$ and 
$\widetilde{\C}^{i,i-1}_{\lambda,\nu}$ $(1 \leq i \leq n)$
are nonzero differential operators of order $\nu-\lambda$ for any $\lambda, \nu \in\C$
with $\nu-\lambda \in \N$. 

The differential operators
$
\widetilde{\C}^{i,i+1}_{\lambda,\nu}$ and 
$
\widetilde{\C}^{i,i-2}_{\lambda,\nu}$ are 
defined only for special parameters $(\lambda,\nu)$ as follows.
\begin{align*}
\widetilde{\C}^{i,i+1}_{\lambda,i+1}&:=
\begin{cases}
 \mathrm{Rest}_{x_n=0} \circ d_{\R^n}&\mathrm{for}\;1\leq i\leq n-2, \lambda=i,\\
d_{\R^{n-1}} \circ \widetilde{\C}_{\lambda,0}&\mathrm{for}\; i=0, \lambda \in -\N,
 \end{cases}\\
\widetilde{\C}^{i,i-2}_{\lambda,n-i+1}&:=
\begin{cases}
\mathrm{Rest}_{x_n=0} 
\circ \iota_{\frac{\partial}{\partial x_n}} d^*_{\R^n}
&\mathrm{for}\; 2\leq i\leq n, \lambda=n-i,\\
-d^*_{\R^{n-1}}\circ \C^{n,n-1}_{\lambda,0}&
\mathrm{for}\; i=n, \lambda \in -\N.
 \end{cases}
\end{align*}
Then we give a complete solution to Problem \ref{prob:2} for 
the model space $(X,Y)=(S^n, S^{n-1})$ in terms of the flat picture of 
principal series representations as follows:
\begin{thm}\label{thm:B}
Suppose a 6-tuple $(i,j,\lambda,\nu,\delta,\eps)$ satisfies one of the equivalent conditions in Theorem \ref{thm:A}. 
Then the 
operators $\widetilde\C^{i,j}_{\lambda,\nu}: C^\infty(\R^n)\otimes\bigwedge^i(\C^n)\To C^\infty(\R^{n-1})\otimes\bigwedge^j(\C^{n-1})$ extend to differential SBOs $I(i,\lambda)_\delta\To J(j,\nu)_\eps$, to be denoted by the same letters. Conversely, any
differential SBO
from $I(i,\lambda)_\delta$ to $J(j,\nu)_\eps$
is proportional to the 
following differential operators:
$
\widetilde{\C}^{i,i}_{\lambda,\nu}$ in Case 1,
$
\widetilde{\C}^{i,i-1}_{\lambda,\nu}$ in Case 2,
$
\widetilde{\C}^{i,i+1}_{i,i+1}$ in Case 3,
$
\widetilde{\C}^{0,1}_{\lambda,1}$ in Case 3$'$,
$
\widetilde{\C}^{i,i-2}_{n-i,n-i+1}$ in Case 4, and
$
\widetilde{\C}^{n,n-2}_{\lambda,1}$ in Case 4$'$.               
\end{thm}

%%%%%%%%%%%%%%%%%%%%%%%%%%%%%%%%%%%%%%%%%%%%%%%
\section{Matrix-valued factorization identities}

 Suppose that
$T_X:I(i,\lambda')_\delta\to I(i,\lambda)_\delta$ or $T_Y:J(j,\nu)_\eps\to J(j,\nu')_\eps$ are $G$- or
$G'$-intertwining operators, respectively. Then
the composition $T_Y\circ D_{X\to Y}$ or $D_{X\to Y}\circ T_X$ of a
symmetry breaking operator $D_{X\to Y}:I(i,\lambda)_\delta\to J(j,\nu)_\eps$
gives another symmetry breaking operator:
$$
\xymatrix{
I(i,\lambda)_\delta\ar[rr]^{D_{X\to Y}}
\ar@{-->}[drr]&&J(j,\nu)_\eps\ar[d]_{T_Y}\\
I(i,\lambda')_\delta\ar[u]^{T_X}\ar@{-->}[urr]&&J(j,\nu')_\eps
}
$$

The multiplicity-free property (see Theorem \ref{thm:A} (ii)) assures the existence of matrix-valued factorization identities
for differential SBOs, namely, $D_{X\to Y}\circ T_X$ must be a scalar multiple of $\widetilde\C_{\lambda',\nu}^{i, j}$, and $T_Y\circ D_{X\to Y}$ must be a scalar multiple of $\widetilde\C_{\lambda,\nu'}^{i, j}$.
We shall determine these constants explicitly when $T_X$ or $T_Y$ are
 Branson's conformally covariant operators \cite{Branson} defined below.
Let $0\leq i\leq n$. For $\ell\in\N_+$, we set
\begin{equation}\label{eqn:T2li}
\mathcal T_{2\ell}^{(i)}:=((\frac n2-i-\ell)d_{\R^{n}}d^*_{\R^{n}}+(\frac n2-i+\ell)d^*_{\R^{n}}d_{\R^{n}})\Delta^{\ell-1}_{\R^n}
=(-2\ell\,d_{\R^{n}}d^*_{\R^{n}}
-(\frac{n}{2}-i+\ell)\Delta_{\R^n})\Delta_{\R^n}^{\ell-1}\nonumber.
\end{equation}

Then the differential operator
$
\mathcal T_{2\ell}^{(i)}:\mathcal E^i(\R^n)\To\mathcal E^i(\R^n)
$
induces a nonzero $O(n+1,1)$-intertwining operator, to be denoted by the same letter $\mathcal T_{2\ell}^{(i)}$, from $I\left(i,\frac n2-\ell\right)_\delta$ to $I\left(i,\frac n2+\ell\right)_\delta$, for $\delta\in\Z/2\Z$.
Similarly, we define a $G'$-intertwining operator
$
\mathcal {T'}_{2\ell}^{(j)}:J\left(j,\frac{n-1}2-\ell\right)_\eps\To J\left(j,\frac{n-1}2+\ell\right)_\eps
$ for $0\leq j\leq n-1$ and $\eps\in\Z/2\Z$
 as the lift of the differential operator $\mathcal {T'}_{2\ell}^{(j)}:\mathcal E^j(\R^{n-1})\To
 \mathcal E^j(\R^{n-1})$ which is given by
\begin{equation*}
\mathcal {T'}_{2\ell}^{(j)}=
((\frac{ n-1}2-j-\ell)d_{\R^{n-1}}d_{\R^{n-1}}^*+(\frac{ n-1}2-j+\ell)d_{\R^{n-1}}^*d_{\R^{n-1}})\Delta^{\ell-1}_{\R^{n-1}}.
\end{equation*}

Consider the following diagrams for $j=i$ and $j=i-1$:
\begin{eqnarray*}
\xymatrix{
I\left(i,\frac n2-\ell\right)_\delta\ar[d]_{\mathcal T_{2\ell}^{(i)}}
\ar[drr]^{{\widetilde \C}^{i,j}_{\frac n2-\ell,\frac n2+a+\ell}}
&& \\
I\left(i,\frac n2+\ell\right)_\delta\ar[rr]_{{\widetilde \C}^{i,j}_{\frac n2+\ell,\frac n2+a+\ell}}&&
J\left(j,\frac n2+a+\ell\right)_\eps,
}
&&
\xymatrix{
I\left(i,\frac {n-1}2-a-\ell\right)_\delta\ar[rr]^{{\widetilde \C}^{i, j}_{\frac{n-1}2-a-\ell,\frac{n-1}2-\ell}}
\ar[drr]_{{\widetilde \C}^{i, j}_{\frac{n-1}{2}-a-\ell,\frac{n-1}2+\ell}}
&&
J\left(j,\frac{n-1}2-\ell\right)_\eps\ar[d]^{{\mathcal T'}^{(j)}_{2\ell}}\\
&&J\left(j,\frac{n-1}2+\ell\right)_\eps,
}
\end{eqnarray*}
where parameters $\delta$ and $\eps \in \Z/2\Z$ are chosen according to Theorem \ref{thm:A} (iii).
In what follows, we put 
$$
p_\pm=
\begin{cases}
i\pm\ell-\frac n2&\mathrm{if}\; a\neq0\\
\pm2&\mathrm{if}\; a=0
\end{cases},
\quad
q=
\begin{cases}
i+\ell-\frac{n-1}2&\mathrm{if}\;i\neq0, a\neq0\\
-2&\mathrm{if}\; i\neq0, a=0\\
-\left(\ell+\frac{n-1}2\right)&\mathrm{if}\; i=0
\end{cases},
\quad
r=
\begin{cases}
i-\ell-\frac{n+1}2&\mathrm{if}\;i\neq n, a\neq0\\
2&\mathrm{if}\; i\neq n, a=0\\
-\left(\ell+\frac{n+1}2\right)&\mathrm{if}\; i=n
\end{cases},
$$
\begin{equation*}
K_{\ell,a}:=\prod_{k=1}^\ell
\left(\left[\frac a2\right]+k\right).
\end{equation*}

Then the factorization identities for differential SBOs $\widetilde\C^{i,j}_{\lambda,\nu}$ for $j\in\{i-1,i\}$ and
Branson's conformally covariant operators $\mathcal T_{2\ell}^{(i)}$ or $\mathcal{T'}_{2\ell}^{(j)}$ are given as follows.
\begin{thm}\label{thm:factor1} Suppose $0\leq i\leq n-1, a\in\N$ and $\ell\in\N_+$.
Then
\begin{eqnarray*}
&(1)&{\widetilde \C}^{i,i}_{\frac n2+\ell,a+\ell+\frac n2}
\circ 
\mathcal {T}_{2\ell}^{(i)}=p_-K_{\ell,a}{\widetilde \C}^{i, i}_{\frac n2-\ell,a+\ell+\frac n2}.\\
&(2)& \mathcal{T'}_{2\ell}^{(i)}\circ{\widetilde \C}^{i,i}_{\frac{n-1}2-a-\ell,\frac{n-1}2-\ell}=
qK_{\ell,a}
{\widetilde \C}^{i,i}_{\frac{n-1}{2}-a-\ell,\frac{n-1}2+\ell}.
\end{eqnarray*}
\end{thm}

\begin{thm}\label{thm:factor2} 
Suppose $1\leq i\leq n$, $a\in\N$ and $\ell\in\N_+$.  Then
\begin{eqnarray*}
&(1)&{\widetilde \C}^{i,i-1}_{\frac n2+\ell,a+\ell+\frac n2}\circ 
\mathcal {T}_{2\ell}^{(i)}=p_+K_{\ell,a}{\widetilde \C}^{i,i-1}_{\frac n2-\ell,a+\ell+\frac n2}.\\
&(2)&\mathcal {T'}_{2\ell}^{(i-1)}\circ 
{\widetilde \C}^{i,i-1}_{\frac{n-1}2-a-\ell,\frac{n-1}2-\ell}=
rK_{\ell,a}{\widetilde \C}^{i,i-1}_{\frac{n-1}{2}-a-\ell,\frac{n-1}2+\ell}.\\
\end{eqnarray*}
\end{thm}
In the case where $i=0$, ${\widetilde \C}^{i, i}_{\lambda,\nu}$ is
a scalar-valued operator, and the corresponding factorization identities in Theorem \ref{thm:factor1} were studied in \cite{Juhl,KS13,KOSS15}.

The main results are proved by using the F-method \cite{K14,KP1,KOSS15}. 
Details will appear elsewhere.

%%%%%%%%%%%%%%%%%%%%%%%%%%%%%%%%%%%%%%%%%%%%%%%

\vskip10pt

\footnotesize{ \noindent  Addresses:
T. Kobayashi.  Kavli IPMU (WPI)
and Graduate School of Mathematical Sciences, 
The University of Tokyo, 3-8-1 Komaba, Meguro, Tokyo, 153-8914 Japan;
\texttt{toshi@ms.u-tokyo.ac.jp}.\vskip5pt

\noindent T. Kubo. Graduate School of Mathematical Sciences, The University of
 Tokyo,
3-8-1 Komaba, Meguro, Tokyo, 153-8914 Japan; \texttt{{
 toskubo@ms.u-tokyo.ac.jp}}.\vskip5pt

\noindent M. Pevzner.
Laboratoire de Math\'ematiques de Reims, Universit\'e
de Reims-Champagne-Ardenne, FR 3399 CNRS, F-51687, Reims, France; \texttt{{
 pevzner@univ-reims.fr.}}

\begin{thebibliography}{99} 

\bibitem{Branson} T.~P. Branson, 
Conformally covariant equations on differential forms,
\emph{Comm. Part. Diff. Eq.} \textbf{7}, (1982),
\href{http://dx.doi.org/10.1080/03605308208820228}
{pp. 393--431}.

\bibitem{FT}
E.~S. Fradkin, A.~A. Tseytlin, 
Asymptotic freedom in extended conformal supergravities,
\emph {Phys. Lett. B} \textbf{110}, (1982)
\href{http://dx.doi.org/10.1016/0370-2693(82)91018-8}
{pp. 117--122}.

\bibitem{GJMS}
C.~R. Graham, R. Jenne, L.~J. Mason, G.~A.~J. Sparling, 
Conformally invariant powers of the Laplacian. I. Existence. 
\emph{J. London Math. Soc.} (2) \textbf{46} (1992), 
\href{http://dx.doi.org/10.1112/jlms/s2-46.3.557}
{pp. 557--565}.

\bibitem{Juhl}A. Juhl, \emph{ Families of conformally covariant differential operators, $Q$-curvature and holography.} Progr. Math., 
\href{http://link.springer.com/book/10.1007/978-3-7643-9900-9/page/1}
{\bf{275}}. Birkh\"auser, Basel, 2009. 

\bibitem{K14} T.~Kobayashi,
F-method for symmetry breaking operators, 
\emph{
Diff. Geometry and its Appl.}
{\textbf{33}}, (2014), 
\href{http://dx.doi.org/10.1016/j.difgeo.2013.10.003}
{pp. 272--289}.

\bibitem{KP1}
T.~Kobayashi, M.~Pevzner,
Differential symmetry breaking operators. I. General theory and
F-method, 
\emph{Selecta. Math. (N.S.)}, \textbf{22}, (2016), 
\href{http://dx.doi.org/10.1007/s00029-015-0207-9}
{pp. 801--845}.

\bibitem{KP2}
T.~Kobayashi, M.~Pevzner, Differential symmetry breaking operators.
II. Rankin--Cohen operators for symmetric pairs, 
\emph{Selecta. Math. (N.S.)}, \textbf{22}, (2016), 
\href{http://dx.doi.org/10.1007/s00029-015-0208-8}
{pp. 847--911}.

\bibitem{KS13}
T.~Kobayashi, B.~Speh,
 Symmetry Breaking for Representations of Rank One Orthogonal Groups,
\href{http://dx.doi.org/10.1090/memo/1126}{Memoirs of American Mathematical Society, vol. \textbf{238}, 2015.} 118 pp.  
{ISBN: 978-1-4704-1922-6}.
 
 \bibitem{KOSS15}
T.~Kobayashi, B.~\O rsted, P. Somberg, and V. Sou\v cek, Branching laws for 
Verma modules and applications in parabolic geometry. I. \emph{Adv. Math}., \textbf{285},
 (2015), 
\href{http://dx.doi.org/10.1016/j.aim.2015.08.020}
{pp. 1796--1852}.

\bibitem{KW14}
B.~Kostant, N. R.~Wallach, Action of the conformal group on steady state solutions to Maxwell's equations and background radiation. Symmetry: representation theory and its applications, 
\href{http://dx.doi.org/10.1007/978-1-4939-1590-3_14}
{pp. 385--418}, 
\emph{Progr. Math.},
{\textbf{257}}, Birkh\"auser/Springer, New York, 2014. 
 
\bibitem{P08}
 S. Paneitz, A quartic conformally covariant differential operator 
 for arbitrary pseudo-Riemannian manifolds,
\emph{SIGMA Symmetry Integrability Geom. Methods Appl.}
\href{http://www.emis.de/journals/SIGMA/2008/}
{\textbf{4} (2008)}, paper 036, 3 pp.

\end{thebibliography}
\end{document}